\documentclass[11pt]{amsart}
\usepackage{a4wide}
\usepackage{amsmath,amsfonts,amssymb}
\usepackage{amsthm}
\usepackage{url}
\usepackage{color}
\usepackage{tikz, tkz-graph}

\newtheorem*{theorem*}{Theorem}
\newtheorem{theorem}{Theorem}
\newtheorem{lemma}[theorem]{Lemma}
\newtheorem{corollary}[theorem]{Corollary}
\newtheorem{proposition}[theorem]{Proposition}

\theoremstyle{remark}

\newcommand{\eps}{\varepsilon}

\begin{document}
\title{On the maximum mean subtree order of trees}

\author[S. Cambie]{Stijn Cambie}
\address{S. Cambie\\
Department of Mathematics\\
Radboud University Nijmegen\\
Postbus 9010\\
6500 GL Nijmegen\\
The Netherlands
} 
\email{s.cambie@math.ru.nl}
\thanks{The first author is supported by a Vidi grant (639.032.614) of the Netherlands Organisation for Scientific Research (NWO)}

\author[S. Wagner]{Stephan Wagner}
\address{S. Wagner\\
Department of Mathematics\\
Uppsala Universitet\\
Box 480\\
751 06 Uppsala\\
Sweden\\
\and
Department of Mathematical Sciences\\ 
Stellenbosch University\\
Private Bag X1\\
Matieland 7602, South Africa\\
}
\email{stephan.wagner@math.uu.se,swagner@sun.ac.za}
\thanks{The second author is supported by the Knut and Alice Wallenberg Foundation.}

\author[H. Wang]{Hua Wang}
\address{Hua Wang\\
Department of Mathematical Sciences \\
Georgia Southern University \\
Statesboro, GA 30460, USA
}
\email{hwang@georgiasouthern.edu}

\begin{abstract}
A subtree of a tree is any induced subgraph that is again a tree (i.e., connected). The mean subtree order of a tree is the average number of vertices of its subtrees. This invariant was first analyzed in the 1980s by Jamison. An intriguing open question raised by Jamison asks whether the maximum of the mean subtree order, given the order of the tree, is always attained by some caterpillar. While we do not completely resolve this conjecture, we find some evidence in its favor by proving different features of trees that attain the maximum. For example, we show that the diameter of a tree of order $n$ with maximum mean subtree order must be very close to $n$. Moreover, we show that the maximum mean subtree order is equal to $n - 2\log_2 n + O(1)$. For the local mean subtree order, which is the average order of all subtrees containing a fixed vertex, we can be even more precise: we show that its maximum is always attained by a broom and that it is equal to $n - \log_2 n + O(1)$.
\end{abstract}

\maketitle

\section{Introduction}

A subtree of a tree $T$ is any induced subgraph that is connected and thus again a tree. In this paper, we will be concerned with the average number of vertices in a subtree (averaged over all subtrees), which is known as the \emph{mean subtree order} of $T$ and denoted $\mu_T$. A normalized version of the mean subtree order, called the \emph{subtree density}, is obtained by dividing by the number of vertices in $T$: if $T$ has $n$ vertices, then
$$D_T = \frac{\mu_T}{n}.$$
This quantity clearly always lies between $0$ and $1$. The concepts of mean subtree order and subtree density were introduced to the literature by Jamison in the 1980s \cite{jamison1983, jamison1984}. Both papers contain a number of interesting questions and conjectures, many of which were only resolved very recently \cite{haslegrave2014, mol2019, vince2010, wagner2014, wagner2016}.

The biggest open problem concerning mean subtree order and subtree density is certainly the natural question: which trees yield the maximum for a given number of vertices? The (highly nontrivial) fact that the minimum is always attained by the path was already  proven by Jamison in his first paper \cite{jamison1983}. Initially, for a small number of vertices, the tree with greatest subtree density is a star. The first exception occurs when the number of vertices $n$ is equal to $9$: here, a double-star obtained by connecting the centers of two stars with four and five vertices respectively has mean subtree order $\frac{779}{159} \approx 4.89937$ and subtree density $\frac{779}{1431} \approx 0.54437$ respectively, compared to the star with mean subtree order $\frac{161}{33} \approx 4.87879$ and density $\frac{161}{297} \approx 0.54209$ respectively.

Based on further computational evidence, Jamison put forward the conjecture that the maximum mean subtree order is attained by a caterpillar (i.e., a tree that becomes a path when all leaves are removed) for every possible number of vertices. This conjecture has been open ever since, and rather little progress has been made. It is not difficult to show that the maximum subtree density approaches $1$ as the number of vertices tends to infinity. There are many possible constructions that yield this limit, the simplest perhaps being ``double-brooms'' (or ``batons''), which consist of a long path with a suitable number of leaves attached to both ends.

Mol and Oellermann \cite{mol2019} considered this construction in greater detail and found that the optimal choice of double-broom for a given number of vertices is essentially to attach (approximately) $2\log_2 n$ leaves at each end of a path of length (approximately) $n - 4\log_2 n$. This immediately yields a lower bound for the maximum of the mean subtree order, which is asymptotically
\begin{equation}\label{eq:moloellerman}
n - 2\log_2 n + O(1),
\end{equation}
as one finds by a relatively straightforward calculation; see also Corollary~\ref{cor:boundsmaxmu} and Theorem~\ref{thr:optimaltreestructure} below for a more precise estimate.

One of the goals of this paper is to show that double-brooms are indeed ``close to optimal''. Specifically, we prove

\begin{theorem}\label{thm:max_asymp}
The maximum of the mean subtree order over all trees with $n$ vertices is $n - 2 \log_2 n + O(1)$.
\end{theorem}

We will see, however, that double-brooms do not attain the maximum for sufficiently large $n$. See Section~\ref{sec:opt_cat} for more details.

Our proof of Theorem~\ref{thm:max_asymp} is, as many other results on the mean subtree order, based on a ``local'' version. Define the local mean subtree order $\mu_T(v)$ at a vertex $v$ of a tree $T$  to be the average number of vertices in a subtree that contains $v$. More generally, one can consider the average of the order of all subtrees containing a specific set $A$ of vertices, denoted $\mu_T(A)$ (so that $\mu_T(v) = \mu_T(\{v\})$ and $\mu_T = \mu_T(\emptyset)$). The following useful monotonicity property was already proven by Jamison \cite{jamison1983}:

\begin{theorem}[{\cite[Theorem 4.5]{jamison1983}}]\label{thm:local_monotonicity}
We have $\mu_T(A) \leq \mu_T(B)$ whenever $A$ is a subset of $B$. Equality holds if and only if the smallest subtree containing all of $A$ is the same as the smallest subtree containing all of $B$.
\end{theorem}

As an important special case \cite[Theorem 3.9]{jamison1983}, 
\begin{equation}\label{eq:localglobal}
\mu_T(v) \geq \mu_T
\end{equation}
holds for every vertex $v$ of $T$ (even with strict inequality unless $T$ only has one vertex). Therefore, any upper bound on local mean subtree orders immediately yields an upper bound on the (global) mean subtree order. The local mean subtree order is often easier to deal with, and in fact we will be able to resolve the local analogue of Jamison's caterpillar conjecture. We even obtain the following stronger result:

\begin{theorem}\label{thm:local_broom}
If the local mean subtree order $\mu_T(r)$ attains its maximum among all choices of an $n$-vertex tree $T$ and a vertex $r$ of the tree, then $T$ has to be a broom, i.e., a tree consisting of a path and leaves attached to one end of the path, and $r$ has to be the other end of the path.
\end{theorem}

The formal proof of this result will be presented in the next section. Let us briefly give an intuitive explanation why brooms are strong candidates for the maximum local mean subtree order. The leaves at one end create a large number of subtrees that all have to contain the end of the path that these leaves are attached to. Since we are only counting subtrees that also contain the other end (the vertex $r$), most of the subtrees under consideration contain the entire path, thus providing a large contribution to the local mean subtree order. One also finds (a precise discussion is given later) that the optimal choice for the number of leaves is about $2 \log_2 n$. The same reasoning also explains why the double-brooms that were mentioned earlier are close to optimal.

Building on Theorem~\ref{thm:local_broom}, Theorem~\ref{thm:max_asymp} will be proven in Section~\ref{sec:maximum}. In Section~\ref{sec:diameter}, we add further evidence in favor of Jamison's caterpillar conjecture by proving that an optimal tree (a tree that attains the maximum mean subtree order) with $n$ vertices, or even a tree that comes close to the maximum, must have a diameter that is very close to $n$. Moreover, we will be able to bound the number of subtrees in an optimal tree both from above and below, showing that it is necessarily of order $\Theta(n^4)$.

Finally, we look closer at caterpillars. We show that with a suitable choice of caterpillar, one can obtain a greater mean subtree order than with a double-broom when the number of vertices is sufficiently large, even though the improvement is modest (of order $O(1)$). The structure of the caterpillars that achieve this feat is somewhat surprising---see Section~\ref{sec:opt_cat} for details.

\section{Extremal local mean subtree order}\label{sec:local}

This section is devoted to the proof of Theorem~\ref{thm:local_broom}. While the calculations are somewhat technical, the main idea is rather straightforward: we prove that a tree that is not a broom can always be improved by replacing a subtree that is the union of two brooms by a single one (keeping the total number of vertices the same) in such a way that the local mean subtree order increases.

\begin{figure}[h]
	\centering
		\begin{tikzpicture}[line cap=round,line join=round,>=triangle 45,x=1.0cm,y=1.0cm]
		\clip(-1,-3) rectangle (3,3);
		\draw [line width=1.pt,dotted] (1.,2.)-- (2.,0.);
		\draw [line width=1.pt,dotted] (1.,2.)-- (-0.5,-1.)-- (1.25,-1.)--(2.,0.);
		\draw [line width=1.pt] (2.,0.)-- (2.,-2.);
		\draw [line width=1.pt] (2.,-2.)-- (1.6,-2.6);
		\draw [line width=1.pt] (2.,-2.)-- (2.4,-2.6);
		\draw [line width=1.pt,dotted] (1.6,-2.6)-- (2.4,-2.6);
		\draw (1,1) node {$T'$};
		\draw (2.,-1) node[anchor= west] {$a$};
		\draw (2,-2.6) node[anchor=north] {$b$};
		\draw (1,2) node[anchor=south] {$r$};
		\draw (2.,0) node[anchor=west] {$v$};
		\begin{scriptsize}
		\draw [fill=black] (2.,0.) circle (2.5pt);
		\draw [fill=black] (1.,2.) circle (2.5pt);
		\draw [fill=black] (2.,-2.) circle (2.5pt);
		\draw [fill=black] (1.6,-2.6) circle (2.5pt);
		\draw [fill=black] (2.4,-2.6) circle (2.5pt);
		\end{scriptsize}
		\end{tikzpicture}
		\quad
		\begin{tikzpicture}[line cap=round,line join=round,>=triangle 45,x=1.0cm,y=1.0cm]
		\clip(-1,-3) rectangle (4,3);
		\draw [line width=1.pt,dotted] (1.,2.)-- (2.,0.);
		\draw [line width=1.pt,dotted] (1.,2.)-- (-0.5,-1.)-- (1.25,-1.)--(2.,0.);
		\draw [line width=1.pt] (2.,0.)-- (1.,-2.);
		\draw [line width=1.pt] (1.,-2.)-- (0.6,-2.6);
		\draw [line width=1.pt] (1.,-2.)-- (1.4,-2.6);
		\draw [line width=1.pt,dotted] (0.6,-2.6)-- (1.4,-2.6);
		\draw [line width=1.pt] (2.,0.)-- (3.,-2.);
		\draw [line width=1.pt] (3.,-2.)-- (2.6,-2.6);
		\draw [line width=1.pt] (3.,-2.)-- (3.4,-2.6);
		\draw [line width=1.pt,dotted] (2.6,-2.6)-- (3.4,-2.6);
		\draw (1,1) node {$T'$};
		\draw (1.5,-1) node[anchor= west] {$a$};		
		\draw (1,-2.6) node[anchor=north] {$b$};
		\draw (2.5,-1) node[anchor= west] {$c$};		
		\draw (3,-2.6) node[anchor=north] {$d$};
		\draw (1,2) node[anchor=south] {$r$};
		\draw (2.,0) node[anchor=west] {$v$};
		\begin{scriptsize}
		\draw [fill=black] (2.,0.) circle (2.5pt);
		\draw [fill=black] (1.,2.) circle (2.5pt);
		\draw [fill=black] (1.,-2.) circle (2.5pt);
		\draw [fill=black] (0.6,-2.6) circle (2.5pt);
		\draw [fill=black] (1.4,-2.6) circle (2.5pt);
		\draw [fill=black] (3.,-2.) circle (2.5pt);
		\draw [fill=black] (2.6,-2.6) circle (2.5pt);
		\draw [fill=black] (3.4,-2.6) circle (2.5pt);
		\end{scriptsize}
		\end{tikzpicture}
	\caption{Abstract configurations $T_1$ and $T_2$}
	\label{fig:configurations}
\end{figure}
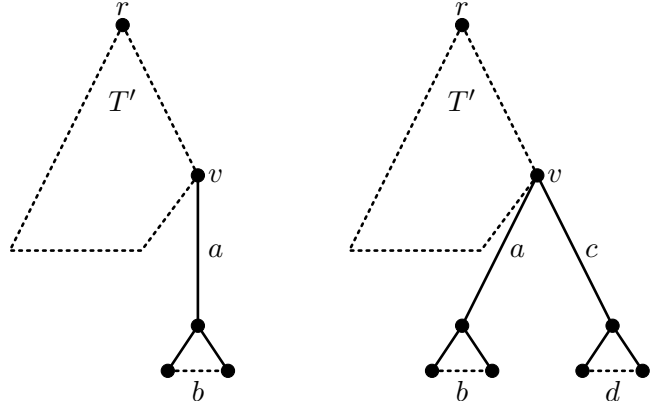

To this end, let us define two general configurations $T_1$ and $T_2$.
The tree $T_1$ consists of a rooted tree $T'$ with root $r$, with a broom attached to a vertex $v$ (of $T'$) whose length is $a\geq0$ and which has $b \geq1$ leaves. Note that $a=0$ is allowed, in which case the broom becomes a star.
Similarly, $T_2$ contains the tree $T'$ with root $r$, with two brooms attached at the same vertex $v$. These brooms have paths of length $a$ and $c$ respectively and $b$ and $d$ leaves respectively. Here we assume $a \geq0$ and $b,c,d \geq1$, since for $a=c=0$ both brooms would degenerate to stars and could be regarded as a single star.
Let us remark that $v$ is not necessarily a leaf of $T'$ in this setup.
The two constructions are presented in Figure~\ref{fig:configurations}. In both cases, we will mainly be interested in the local mean subtree order at the root.

So let $\ell$ be the number of subtrees of $T'$ that contain $r$ and $v$, and let $s$ be their total order (number of vertices).
Likewise, let $m$ be the number of subtrees of $T'$ that contain $r$, but not $v$, and let $t$ be their total order. Note that $\mu_{T'}(\{r,v\})=\frac{s}{\ell}$ and $\mu_{T'}(\{r\})=\frac{s+t}{\ell+m}$ by definition, and recall from Theorem~\ref{thm:local_monotonicity} that $\mu_{T'}(\{r,v\}) \geq \mu_{T'}(\{r\})$ where equality only holds when $v = r$.
This is equivalent to
$$\frac{s}{\ell} \geq\frac{s+t}{\ell+m} \Leftrightarrow sm \geq t \ell \Leftrightarrow \frac{s}{\ell} \geq \frac t m .$$

We can now express the local mean subtree order of $T_1$ at $r$ in terms of the variables $a,b,\ell,m,s,t$. Note that there are $m$ subtrees of $T_1$ containing $r$, but not $v$, and $\ell(a+2^b)$ subtrees of $T_1$ containing both $r$ and $v$: any subtree of $T'$ that contains both vertices can be combined with any of the $a+2^b$ subtrees of the attached broom that contains $v$ (here, $a$ is the number of proper subpaths of the path of length $a$, and $2^b$ is the number of subtrees containing the entire path of length $a$ and a subset of the $b$ leaves). The total order of the former is $t$ by definition, the total order of the latter is $s(a+2^b) + \frac{\ell}{2}(a^2-a+(2a+b)2^b)$ (the first term being the contribution from $T'$, the latter the contribution from the broom). It follows that
\begin{equation}\label{eq:meanorderT1}
\mu_{T_1}(r) = \frac{t+s(a+2^b)+\frac {\ell}2 \left( a^2-a+(2a+b)2^b \right) }{m+\ell(a+2^b)}.
\end{equation}
The number of variables can be reduced using the abbreviations $k= \frac{m}{\ell}$ and $C=2 \frac{sm-t\ell}{\ell^2}$:
$$\mu_{T_1}(r) =\frac s{\ell} + \frac 12 \frac{  a^2-a+(2a+b)2^b -C }{k+(a+2^b)}.$$
Note here that $C \geq 0$ by the earlier observation that $\frac{s}{\ell} \geq \frac t m$. Using analogous reasoning, we obtain 
\begin{align*}
\mu_{T_2}(r)
&= \frac{t+s(a+2^b)(c+2^d)+\frac {\ell}2 \left( (a+2^b)(c^2-c+(2c+d)2^d) +(c+2^d)(a^2-a+(2a+b)2^b) \right) }{m+\ell(a+2^b)(c+2^d)}\\
&=\frac s{\ell} + \frac 12 \frac{  (a+2^b)(c^2-c+(2c+d)2^d) +(c+2^d)(a^2-a+(2a+b)2^b) -C }{k+(a+2^b)(c+2^d)}.\end{align*}

Since $\frac s{\ell}$ is determined by the structure of $T'$ alone, we focus on the remaining parts of $\mu_{T_1}(r)$ and $\mu_{T_2}(r)$ and write these as functions $f_1$ and $f_2$ of the six variables $k,C,a,b,c,d$, i.e.,
\begin{align*}
f_1(k,C,a,b)&= \frac{  a^2-a+(2a+b)2^b -C }{k+(a+2^b)}\mbox{ and} \\
	f_2(k,C,a,b,c,d)&=\frac{  (a+2^b)(c^2-c+(2c+d)2^d) +(c+2^d)(a^2-a+(2a+b)2^b) -C }{k+(a+2^b)(c+2^d)}.
\end{align*}

Ultimately, we would like to obtain an inequality of the form $f_1(k,C,\alpha,\beta) \geq f_2(k,C,a,b,c,d)$ for a suitable choice of $\alpha,\beta$ that may depend on $k,C,a,b,c,d$.  We first prove some auxiliary lemmas for this purpose.

\begin{lemma}\label{lem1}
	Let $k, C \geq 0$ be fixed constants, and let $N \geq 2$ be a fixed integer. Suppose that the integers $a$ and $b$ maximize the function $f_1(k,C,a,b)$ under the conditions $a+b=N$, $a \geq 0$ and $b \geq 1$. Then we have $2^b \geq 3a$, except when $k=C=0$ and $N=2$, in which case $a=0, b=2$ and $a=b=1$ both maximize $f_1(k,C,a,b)$.
\end{lemma}

\begin{proof}
	Notice that $a+2^b=N-b+2^b$ is an increasing function of $b$, so  $\frac{  -C }{k+(a+2^b)}$ is increasing in $b$ (or constant if $C=0$), and 
	$ \frac{a+2^b}{k+(a+2^b)}$ is increasing in $b$ as well when $k>0.$
	So if a certain choice of $b$ maximizes the function $f_1(0,0,N-b,b)$, then $f_1(k,C,N-b,b)=\frac{a+2^b}{k+(a+2^b)} f_1(0,0,a,b)+\frac{  -C }{k+(a+2^b)}$ cannot attain its maximum for any smaller value of $b$ when $(k,C) \not = (0,0)$.
	So it suffices to prove the statement in the case that $k=C=0$.
	
	Note that $$f_1(0,0,a,b) = \frac{  a^2-a+(2a+b)2^b }{(a+2^b)} = 2N-b - \frac{ a(N+1)}{a+2^b}.$$
	If $2^b < 3a$, then this is strictly smaller than $2N-1-\frac{N+1}4$. On the other hand, if we choose $b= \lfloor 2 \log_2 N \rfloor $ (and $a = N - b$), then
$$\frac{ a(N+1)}{a+2^b} < \frac{a(N+1)}{a+2^{2\log_2 N-1}} = \frac{2aN+2a}{2a+N^2} < \frac{2N^2+4a}{2a+N^2} = 2$$
and thus
$$2N - b - \frac{ a(N+1)}{a+2^b} > 2N - \lfloor 2 \log_2 N \rfloor - 2.$$
For $N  \geq 47$, we have $2N-1-\frac{N+1}4 \leq  2N - \lfloor 2 \log_2 N \rfloor - 2$, and the conclusion readily follows. For $3 \leq N \leq 46$, the claim is easily checked with a small computer program\footnote{\url{https://github.com/StijnCambie/Jamison}, document Lemma1Check.}.
For $N=2$, it is easily checked by hand that $f_1(0,0,0,2)=f_1(0,0,1,1)=2$ and when $\max\{k,C\}>0$ we have $f_1(k,C,0,2)>f_1(k,C,1,1).$
\end{proof}

\begin{lemma}\label{pos}
	For real numbers $a,c,x,y$ for which $x \geq 3a$, $y \geq 3c$ as well as $a\geq 0$ and $c \geq 1$,
	the expression $xy-xc-ay-(c-1)(a-1)$ is nonnegative.
\end{lemma}
\begin{proof} We have
\begin{align*}
xy-xc-ay-(c-1)(a-1) &= (x-a)(y-c) - 2ac + a + c - 1 \\
&\geq 2a \cdot 2c - 2ac + a + c - 1 \\
&= 2ac + a + c - 1  \geq 0. \qedhere
\end{align*}
\end{proof}

\begin{proposition}\label{2to1}
	Let $k, C \geq 0$ be fixed constants. Let $a,b,c,d$ be integers with $a \geq 0$ and $b,c,d \geq 1$ such that $2^b \geq 3a$ and $2^d \geq 3c$ holds.
Then we have
$$\max\{ f_1(k,C,a+c,b+d) , f_1(k,C,a+c-1,b+d+1)\} > f_2(k,C,a,b,c,d).$$
\end{proposition}

\begin{proof}
We prove that there is a linear combination with nonnegative coefficients of $\Delta_1= f_1(k,C,a+c,b+d)   -f_2(k,C,a,b,c,d)$ and $\Delta_2=f_1(k,C,a+c-1,b+d+1) -f_2(k,C,a,b,c,d)$ which is strictly positive. This implies that either $\Delta_1 > 0$ or $\Delta_2 > 0$, from which the result follows. Let us write $x=2^b$ and $y=2^d$. We take the following coefficients: 
\begin{align*}
\lambda_1 &=\left(xy-xc-ay-(c-1)(a-1)\right)(k+a+c+xy), \\
\lambda_2 &=\left(xc+ay+(c-1)a-c\right)(k+a+c-1+2xy).
\end{align*}
Note that $ \lambda_1 \geq 0$ by Lemma~\ref{pos}, and clearly also $\lambda_2>0.$

The linear combination\footnote{A verification can be found at \url{https://github.com/StijnCambie/Jamison}, document Proposition1Check.} 
\begin{align*}
\lambda_1 \Delta_1 + \lambda_2 \Delta_2=& 
c\left( c-1 \right) \left( xy-1 \right)   (x+a-1)+a \left( (a-1 )(xy-1)+bxy \right) \left(c+    y-1  \right)\\&+cyx \left( (x-1)d-b \right) 
+yad \left( x(c-1)+1 \right) +bxc
\end{align*}
is indeed strictly positive for all integers $a,b,c,d,x,y$ with $a \geq 0$ and $b,c,d \geq 1$ as well as $x \geq 1+b \geq 2$ and $y \geq 2$. It follows that either $\Delta_1 > 0$ or $\Delta_2 > 0$ (or both), completing the proof.
\end{proof}

We are now ready to prove the main result of this section.

\begin{proof}[Proof of Theorem~\ref{thm:local_broom}]
Suppose there is an optimal tree for the local mean subtree order with respect to the vertex $r$ (the root) that is not a broom.
Then there is a vertex $v$, possibly equal to $r$, which is the root of at least two brooms (take the vertex at greatest distance from the root for which the tree consisting of this vertex and all its successors is not a broom). 
Hence the optimal tree can be described in the way we defined the general construction $T_2$.

At the same time, it can also be regarded as a $T_1$ (with the other broom becoming part of $T'$), so Lemma~\ref{lem1} applies: $a,b$ and $c,d$ have to maximize $f_1(k,C,.,.)$ for the corresponding $k$ and $C$ under the fixed sum condition as the tree is assumed to be optimal. We can conclude that $2^b \geq 3a$, unless $k = C = 0$ (which happens only if $v=r$) and $N=2$. In that case, we may assume $a=0$ and $b=2$ without loss of generality as both this choice and $a=b=1$ yield the same local density. So we always have $2^b \geq 3a$, and likewise $2^d \geq 3c$. 
We can also assume $a \geq 0$ and $b,c,d \geq 1$ as mentioned before.

But now we can apply Proposition~\ref{2to1} and conclude that we can replace the two brooms by a single one such that the order is the same and the local mean subtree order increases.
So the original tree was not optimal, which contradicts the initial assumption.
We conclude that an optimal tree has to be a broom as described in the statement of the theorem.
\end{proof}

\begin{corollary}\label{sharplocalresult}
	Among all trees of order $n$, the maximum local mean subtree order is of the form $n- \log_2 n + O(1).$
	More precisely, it is of the form 
$$n- \log_2 n + \frac12 f(2\log_2 n) + o(1),$$
where $f$ is the $1$-periodic function given by $f(x) = x - 2^x$ for $x \in [0,1]$.
\end{corollary}

\begin{proof}
	We know that the maximum local subtree order is attained by a broom. Let $a$ denote the length of its ``handle'' and $b$ the number of leaves, so that $n = a + b + 1$. We can take $\ell = 1$, $m = 0$, $s = 1$ and $t = 0$ in~\eqref{eq:meanorderT1} to see that the average subtree order is $n- \frac 12 \left( b + \frac{an}{a+2^b}\right)$. Hence we need to minimize $b + \frac{an}{a+2^b}$, subject to the condition that $a+b = n-1$.

If we take $b = \lceil 2\log_2 n \rceil$, then $\frac{an}{a+2^b} \leq \frac{an}{2^b} \leq \frac{n^2}{n^2} = 1$. So the minimum of our expression is at most $\lceil 2 \log_2 n \rceil + 1 \leq 2 \log_2 n +2$. Consequently, any choice of $a$ and $b$ where $b > 2\log_2 n + 2$ cannot be optimal. It follows that $a = n - O(\log n)$ for the optimal choice of $a$ and $b$. 

Likewise, we must have
$$\frac{an}{a+2^b} \leq 2\log_2 n + 2$$
for the optimal choice of $a$ and $b$, which implies that $2^b \geq \frac{n^2}{2\log_2 n+2} - O(n)$. From these estimates, we obtain
$$\frac{an}{a+2^b} = \frac{n^2}{2^b} + O \Big( \frac{\log^2 n}{n} \Big).$$
In other words, the minimum of $b + \frac{an}{a+2^b}$, subject to the condition that $a+b = n-1$, differs from the minimum of $b + \frac{n^2}{2^b}$ only by $o(1)$.

The minimum of $b + \frac{n^2}{2^b}$ occurs (by convexity) for a value of $b$ for which we simultaneously have $b + \frac{n^2}{2^b} \leq b-1 + \frac{n^2}{2^{b-1}}$ and $b + \frac{n^2}{2^b} \leq b+1 + \frac{n^2}{2^{b+1}}$.
Those yield $2^b \leq n^2 \leq 2^{b+1}$, i.e., we can take $b = \lfloor 2 \log_2 n \rfloor$. 

We obtain that the maximum local subtree order of a broom (and thus an arbitrary tree) of order $n$ is
$$n- \frac 12 \left( b + \frac{an}{a+2^b}\right) = n - \frac{b}{2} - \frac{n^2}{2^{b+1}} + o(1) = n - \frac{2 \log_2 n - x}{2} - \frac{n^2}{2^{2\log_2 n-x+1}} + o(1),$$
where $x = \{2\log_2 n\}$ is the fractional part of $2\log_2 n$. This can also be written as
$$n - \log_2 n + \frac12 (x - 2^x) + o(1),$$
which completes the proof of our asymptotic formula.
\end{proof}

In the following section, we will use our knowledge on the maximum local mean subtree order to bound the global mean subtree order as well.

\section{The maximum mean subtree order}\label{sec:maximum}

We now make the step from the local to the global mean subtree order. Recall that the global mean subtree order of a tree $T$ is no greater than the local mean subtree order at any vertex of $T$ (inequality~\eqref{eq:localglobal}). This combined with Corollary~\ref{sharplocalresult} shows immediately that
$$\mu_T \leq n - \log_2 n + O(1)$$
for every tree of order $n$. However, in order to match the lower bound~\eqref{eq:moloellerman} due to Mol and Oellermann, we have to refine the argument to prove our main result on the mean subtree order (Theorem~\ref{thm:max_asymp}).

A vertex of a tree is said to be a centroid vertex if none of the components that result when the vertex is removed contains more than half of the vertices. It is well known that every tree has either one or two centroid vertices, and that a centroid vertex minimizes the sum of all distances to the other vertices \cite{zelinka1968}.

\begin{proposition}\label{prop:centroid}
Let $T$ be a tree of order $n$, and let $v$ be a centroid vertex of $T$. Then the local mean subtree order of $T$ at $v$ is at most
$$n - 2\log_2 n + 2 + f(2\log_2 n) + o(1),$$
with $f$ as in Corollary~\ref{sharplocalresult}.
\end{proposition}

\begin{proof}
We can assume that the local mean subtree order of $T$ at $v$ is greatest among all choices of a tree $T$ and a centroid vertex $v$. Let the components of $T - v$ be $T_1,T_2,\ldots,T_k$. Moreover, let $T_i'$ be the tree that results when $v$ is added back to $T_i$ (along with the edge connecting it to its neighbor in $T_i$). Clearly, every subtree of $T$ that contains $v$ induces a subtree in each $T_i'$ that contains $v$, and conversely every subtree of $T$ that contains $v$ can be obtained by merging subtrees of $T_1',T_2',\ldots,T_k'$. It follows easily from this observation that
\begin{equation}\label{eq:mu_add}
\mu_T(v) - 1 = \sum_{i=1}^k \big(\mu_{T_i'}(v) - 1 \big).
\end{equation}
Note also that $\sum_{i=1}^k |T_i| = \sum_{i=1}^k (|T_i'| - 1) = n-1$. We use this representation combined with the upper bound on the local mean subtree order to bound the local mean subtree order at the centroid vertex $v$. Assume without loss of generality that $|T_1| \leq |T_2| \leq \cdots \leq |T_k|$. 

If $k \geq 4$, then $|T_1| \leq |T_2| \leq \frac{n-1}{4}$, thus $|T_1' \cup T_2'| = |T_1| + |T_2| + 1 \leq \frac{n+1}{2}$. By Theorem~\ref{thm:local_broom}, we can replace $T_1' \cup T_2'$ by a single broom of the same order whose local mean subtree order at the root is greater than the local mean subtree order of $T_1' \cup T_2'$. Thus the local mean subtree order increases after this replacement, while vertex $v$ remains a centroid since $T_1$ and $T_2$ together do not contain more than half of the vertices. This contradicts our choice of $T$ and $v$. So we know that $k \leq 3$.

If $k = 2$, then $T_1'$ and $T_2'$ each contain between $\frac{n}{2}$ and $\frac{n}{2}+1$ vertices, and we can apply~\eqref{eq:mu_add} and Corollary~\ref{sharplocalresult} to complete the proof. If $k = 3$, then we note first that $|T_3'| \geq |T_2'| \geq \frac{n}{4} + \frac12$ (since $T_3$ contains at most half of the vertices, and $T_2'$ more than half of the rest). If $|T_1'| \geq \sqrt{n}$, then we can already apply~\eqref{eq:mu_add} and Corollary~\ref{sharplocalresult} to $T_1'$, $T_2'$ and $T_3'$ to show that the local mean subtree order at $v$ is at most
$$n - \frac52 \log_2 n + O(1).$$
In this case, we are done. Otherwise, $T_2'$ and $T_3'$ each contain at least $\frac{n}{2} - \sqrt{n}$ vertices (and at most $\frac{n}{2} + 1$), by the choice of $v$ as a centroid. Applying Corollary ~\ref{sharplocalresult} to those two, we find that the local mean subtree order at $v$ is indeed at most
$$n - 2\log_2 \Big( \frac{n}{2} - \sqrt{n} \Big) + f\Big( 2\log_2 \Big(\frac{n}{2} - \sqrt{n} \Big) \Big) + o(1) = n - 2\log_2 n + 2 + f(2\log_2 n) + o(1).$$
Note here that $\log_2 (\frac{n}{2} - \sqrt{n}) = \log_2 n - 1  +O(n^{-1/2})$, and that $f$ is a continuous $1$-periodic function. This completes the proof.
\end{proof}

Proposition~\ref{prop:centroid} tells us that every tree has a vertex such that the local mean subtree order at that vertex is at most $n - 2 \log_2 n + O(1)$. Specifically, every centroid vertex has this property. Combined with~\eqref{eq:localglobal}, this immediately implies Theorem~\ref{thm:max_asymp}.

It is worth studying the construction of Mol and Oellermann given in \cite{mol2019} a little further: if we take a double-broom consisting of a path of $n-2s$ vertices, with $s$ leaves attached at each end, then there are
$$2^{2s} + 2^{s+1}(n-2s-1) + \binom{n-2s-1}{2} + 2s$$
subtrees with a total of
$$2^{2s}(n-s) + 2^s (n-s)(n-2s-1) + \binom{n-2s}{3} + 2s$$
vertices. The optimal choice of $s$ is $2\log_2 n + O(1)$, in which case we have $2^s = \Theta(n^2)$. The mean subtree order becomes
\begin{align*}
\frac{2^{2s}(n-s) + 2^s (n-s)(n-2s-1) + \binom{n-2s}{3} + 2s}{2^{2s} + 2^{s+1}(n-2s-1) + \binom{n-2s-1}{2} + 2s} &= \frac{2^{2s}(n-s) + 2^s n^2 + O(n^3 \log n)}{2^{2s} + 2^{s+1}n + O(n^2 \log n)} \\
&= n - s - \frac{n^2}{2^s} + O \Big( \frac{\log n}{n} \Big).
\end{align*}

The expression $s + \frac{n^2}{2^s}$ was already analyzed in the proof of Corollary~\ref{sharplocalresult}. Its minimum is obtained when $s = \lfloor 2\log_2 n \rfloor$ and has a value of $2\log_2 n - f(2\log_2 n)$, with the same function $f$ as in that corollary. So we have the following refinement of Theorem~\ref{thm:max_asymp}:

\begin{corollary}\label{cor:boundsmaxmu}
The maximum of the mean subtree order over all trees with $n$ vertices lies between $n - 2 \log_2 n + f(2\log_2 n) + o(1)$ and $n - 2 \log_2 n + 2 + f(2\log_2 n) + o(1)$.
\end{corollary}

We have thus reduced the gap between the upper and lower bound to $2+o(1)$.

\section{Characteristics of optimal trees}\label{sec:diameter}

We call a tree \emph{optimal} if it has the greatest possible mean subtree order among all trees of the same size. In this section, we prove several structural properties of optimal trees. Specifically, we show that the diameter of optimal trees with $n$ vertices is very close to $n$, and that the number of subtrees is of the order of magnitude $\Theta(n^4)$. Another feature of optimal trees was proved by Mol and Oellermann in \cite{mol2019}, namely that the number of leaves cannot be too large:

\begin{theorem}[{\cite[Corollary 4.4]{mol2019}}]
The number of leaves in an optimal tree with $n$ vertices is at most $4 \log_2 n + O(1)$.
\end{theorem}

This estimate is based on a result of Jamison \cite[Lemma 6.1]{jamison1983} stating that the mean subtree order of a tree with $n$ vertices and $\ell$ leaves is at most $n - \ell/2$.

The following notation will be useful for our arguments. For an arbitrary tree $T$, we denote the total number of subtrees by $\sigma(T)$, and their total order (i.e., the total number of vertices in all these subtrees) by $\tau(T)$, so that the mean subtree order is given by
$$\mu_T = \frac{\tau(T)}{\sigma(T)}.$$
For a rooted tree $T$ with root $r$, we let $s(T)$ be the number of subtrees of $T$ that contain the root, and let $t(T)$ be their total order (we suppress the dependence on the root for simplicity). The local subtree order at $r$ can then be expressed as the quotient of these two:
$$\mu_{T}(r) = \frac{t(T)}{s(T)}.$$
We also define the \emph{defect} $\Delta(T)$ as the average number of vertices \emph{not} contained in a randomly chosen subtree that contains the root $r$, i.e.,
$$\Delta(T) = |T| - \mu_{T}(r) = |T| - \frac{t(T)}{s(T)}.$$
Note in particular that the defect of a single-vertex tree is $0$ while the defect of a (rooted) tree with exactly two vertices is $\frac12$. Let $T_1,T_2,\ldots,T_k$ be the 
components resulting when the root is removed. Each of them is endowed with a natural root, namely the unique neighbor of $T$'s root. Let moreover $T_1',T_2',\ldots,T_k'$ be the trees obtained from $T_1,T_2,\ldots,T_k$ by adding back the root of $T$ (and an edge connecting it to the root of the respective component, similar to the proof of Proposition~\ref{prop:centroid}). We have
\begin{equation}\label{eq:s_rec}
s(T) = \prod_{i=1}^k s(T_i'),
\end{equation}
since each subtree that contains the root of $T$ can be decomposed into subtrees of $T_1',T_2',\ldots,T_k'$ in a natural way. Furthermore, when a uniformly random subtree of $T$ that contains the root is chosen, the induced subtrees in $T_1',T_2',\ldots,T_k'$ are independent uniformly random subtrees containing the root in the respective branches. It follows that
\begin{equation}\label{eq:delta_rec}
\Delta(T) = \sum_{i=1}^k \Delta(T_i').
\end{equation}

The following inequality between $\Delta(T)$ and $s(T)$ will be crucial:

\begin{lemma}\label{lem:logsT}
For every rooted tree $T$, we have
$$\Delta(T) \geq \frac12 \log_2 s(T).$$
\end{lemma}

\begin{proof}
We proceed by induction on $|T|$. For a single-vertex tree, the inequality is trivial as both sides are equal to $0$. Next we distinguish two different cases: if there are two or more branches, then we can use the induction hypothesis together with~\eqref{eq:s_rec} and~\eqref{eq:delta_rec}. Otherwise, there is only a single component $T_1$, and we obtain $s(T) = 1 + s(T_1)$, $t(T) = 1+ s(T_1) + t(T_1)$, and consequently
\begin{align*}
\Delta(T) &= |T_1| - \frac{t(T_1)}{1+s(T_1)} = \Delta(T_1) + \frac{t(T_1)}{s(T_1)(1+s(T_1))} \\
&\geq \Delta(T_1) + \frac{1}{1+s(T_1)}.
\end{align*}
On the other hand,
\begin{align*}
\frac12 \log_2 s(T) &= \frac12 \log_2 (1+s(T_1)) = \frac12 \log_2 s(T_1) + \frac12 \log_2 \Big(1 + \frac{1}{s(T_1)} \Big) \\
&\leq \frac12 \log_2 s(T_1) +  \frac{1}{1+s(T_1)}. 
\end{align*}
Thus we can again invoke the induction hypothesis to complete the proof of the inequality.
\end{proof}

As an immediate consequence of Lemma~\ref{lem:logsT}, we already get an upper bound on the number of subtrees of an optimal tree.

\begin{proposition}\label{prop:subtrees_upper}
There is a constant $C_1 > 0$ such that every optimal tree $T$ has at most $C_1 n^4$ subtrees, where $n$ is the number of vertices of $T$.
\end{proposition}

\begin{proof}
It is clearly enough to prove the statement for sufficiently large $n$. Let $T$ be an optimal tree, and suppose first that none of the vertices of $T$ is contained in more than half of the subtrees. Then the mean subtree order is clearly at most
$\frac{n}{2}$, which contradicts Theorem~\ref{thm:max_asymp} (at least for sufficiently large $n$). Thus we can select a vertex $r$ as the root that is contained in more than $\sigma(T)/2$ subtrees. By~\eqref{eq:localglobal} and Lemma~\ref{lem:logsT}, we have
$$\mu_T \leq \mu_T(r) = n - \Delta(T) \leq n - \frac12 \log_2 s(T) \leq n - \frac12 \log_2 \sigma(T) + \frac12.$$
Combining this inequality with Theorem~\ref{thm:max_asymp}, we obtain $\log_2 \sigma(T) \leq 4\log_2 n + O(1)$, which implies the statement.
\end{proof}

Next we show that every optimal tree has a large central part that is contained in most subtrees. Formally, we define the central part $C(T)$ of a tree $T$ to be the set of all vertices that are contained in at least $\frac{1}{1+n^{-1/4}} \sigma(T)$ of all subtrees. The constant $\frac14$ is somewhat arbitrary in this definition and can be replaced by any other number less than $\frac12$. We remark that the definition of the central part is conceptually similar to that of the \emph{subtree core} as defined in \cite{szekely2005}: the subtree core contains those vertices that are contained in the greatest number of subtrees. It can be shown that there are always either one or two vertices in the subtree core; our central part turns out to be much larger for optimal trees.

\begin{lemma}\label{lem:central_part}
Let $T$ be an optimal tree with $n$ vertices, where $n$ is sufficiently large. The vertices of the central part $C(T)$ induce a connected graph, i.e., a subtree of $T$, with at least $n - n^{1/3}$ vertices. Moreover, the subtree induced by $C(T)$ has at most $16$ leaves.
\end{lemma}

\begin{proof}
Let $\sigma_v(T)$ denote the number of subtrees containing a vertex $v$, so that $C(T)$ consists of all vertices for which $\sigma_v(T) \geq \frac{1}{1+n^{-1/4}} \sigma(T)$. It was shown in \cite[Theorem 9.1]{szekely2005} that $\sigma_v(T)$ is unimodal along paths: for any path from one leaf to another, it first increases, then decreases. Hence the minimum of $\sigma_v(T)$ among all vertices on an arbitrary path is attained at one or both ends. Consequently, for any two vertices $v$ and $w$ that belong to $C(T)$, the entire path between $v$ and $w$ is also contained in $C(T)$. This proves that $C(T)$ induces a connected graph. With some minor abuse of notation, we will also write $C(T)$ for the graph induced by $C(T)$. To bound the number of vertices in $C(T)$, we use the following crude bounds:
\begin{itemize}
\item Every vertex outside of $C(T)$ is contained in at most $\frac{1}{1+n^{-1/4}} \sigma(T)$ subtrees.
\item Every vertex in $C(T)$ is contained in at most $\sigma(T)$ subtrees.
\end{itemize}
Thus
$$\tau(T)  \leq (n - |C(T)|) \cdot \frac{1}{1+n^{-1/4}} \sigma(T) + |C(T)| \sigma(T) = \Big( n - \frac{n-|C(T)|}{n^{1/4} + 1} \Big) \sigma(T).$$
So by Theorem~\ref{thm:max_asymp},
$$n - 2\log_2 n + O(1) = \mu(T) = \frac{\tau(T)}{\sigma(T)} \leq n - \frac{n-|C(T)|}{n^{1/4} + 1},$$
which implies that
$$n - |C(T)| \leq 2n^{1/4} \log_2 n + O(n^{1/4}) \leq n^{1/3}$$
for sufficiently large $n$, so that $|C(T)| \geq n - n^{1/3}$. It remains to prove the statement on the number of leaves. 

Consider any leaf $v$ of the subtree $C(T)$ and let $w$ be its unique neighbor in $C(T)$.  When the edge $vw$ is removed from $T$, we obtain two components. The component containing $v$ is called the \emph{branch bundle} of $v$, denoted $B(v)$; the branch bundle $B(v)$ has $v$ as a natural root. 

Let $a$ be the number of subtrees of $T - B(v)$ that contain $w$, and let $b = s(B(v))$ be the number of subtrees of $B(v)$ that contain $v$. Then the total number of subtrees of $T$ that contain $v$ is $(a+1)b$, as every subtree of $B(v)$ containing $v$ is such a subtree itself, and can also be combined with an arbitrary subtree of $T - B(v)$ containing $w$. We observe that $T$ has at least $a+1$ subtrees that do not contain $v$: all $a$ subtrees of $T - B(v)$ that contain $w$, and any single vertex of $T$ other than $v$ and $w$ can also be regarded as such a subtree. Thus we have
$$\sigma(T) \geq (a+1) + (a+1)b = (a+1)(b+1),$$
which implies
$$(a+1)b = (a+1)(b+1) \cdot \frac{b}{b+1} \leq \sigma(T) \cdot \frac{b}{b+1}.$$
By definition of the central part $C(T)$, we must have $(a+1)b \geq \frac{1}{1+n^{-1/4}} \sigma(T)$, so
$$\frac{1}{1+n^{-1/4}} \sigma(T) \leq \sigma(T) \cdot \frac{b}{b+1},$$
which is ultimately equivalent to 
\begin{equation}\label{eq:st_branch}
b =s(B(v)) \geq n^{1/4}. 
\end{equation}
Now suppose that the subtree $C(T)$ has more than 16 leaves. The branch bundles associated with these leaves are all disjoint, since each of them only contains one vertex of $C(T)$. We can choose a subtree containing the root in each of these branch bundles and take their union with the tree $C(T)$ to obtain a subtree of $T$. In view of~\eqref{eq:st_branch}, this gives us at least $(n^{1/4})^{17} = n^{4.25}$ different subtrees of $T$, which contradicts Proposition~\ref{prop:subtrees_upper} for sufficiently large $n$.
\end{proof}

Now we can obtain a matching lower bound for Proposition~\ref{prop:subtrees_upper}, showing that an optimal tree with $n$ vertices has $\Theta(n^4)$ subtrees.

\begin{proposition}\label{prop:subtrees_lower}
There is a constant $C_2 > 0$ such that every optimal tree $T$ has at least $C_2 n^4$ subtrees, where $n$ is the number of vertices of $T$.
\end{proposition}

\begin{proof}
We build our proof on the structural properties established in Lemma~\ref{lem:central_part}. Once again, it suffices to prove the statement for sufficiently large $n$. Consider two families of subtrees of an optimal tree $T$ with $n$ vertices:
\begin{itemize}
\item Let $\mathcal{F}_1$ be the family of subtrees that contain all leaves of the tree $C(T)$, and thus all of $C(T)$. Note that each of these leaves is, by definition, contained in all but at most
$$\sigma(T) -  \frac{1}{1+n^{-1/4}} \sigma(T)  = \frac{\sigma(T)}{n^{1/4}+1}$$
of $T$'s subtrees. Thus the number of subtrees of $T$ that do not contain one of these leaves is at most
$$\frac{16\sigma(T)}{n^{1/4}+1},$$
which in turn means that $\mathcal{F}_1$ contains at least $(1 - 16n^{-1/4})\sigma(T)$ subtrees. If we contract all vertices of $C(T)$ to a single vertex $r$ and call the resulting tree $T'$, then each of the subtrees in $\mathcal{F}_1$ becomes a subtree of $T'$ that contains $r$ (which we take as the root of $T'$), and this correspondence is clearly bijective. The average number of vertices not contained in subtrees of $\mathcal{F}_1$ thus becomes exactly the defect $\Delta(T')$, which we can bound in the following way by Lemma~\ref{lem:logsT}:
\begin{align*}
    \Delta(T') & \geq \frac12 \log_2 s(T') \\ & \geq \frac12 \log_2 \Big( (1 - 16n^{-1/4})\sigma(T) \Big) \\
&= \frac12 \log_2 \sigma(T) - O(n^{-1/4}).
    \end{align*}
In summary, we find that the total number of vertices \emph{not} contained in subtrees that belong to $\mathcal{F}_1$ is at least
\begin{align*} |\mathcal{F}_1| \Delta(T') & \geq (1 - 16n^{-1/4})\sigma(T) \Big( \frac12 \log_2 \sigma(T) - O(n^{-1/4}) \Big) \\
& = \frac{\sigma(T)}{2} \log_2 \sigma(T) \big(1 - O(n^{-1/4}) \big).
\end{align*}
\item In order to define the second family, we consider a centroid vertex $v$ of $T$ (cf. Proposition~\ref{prop:centroid}). Let the components of $T - v$ be $T_1,T_2,\ldots,T_k$, and let $T_1', T_2', \ldots, T_k'$ be obtained from these components by adding back the vertex $v$ in the same way as in Proposition~\ref{prop:centroid}. Recall that none of $T_1,T_2,\ldots,T_k$ can contain more than half of the vertices of $T$. Since we know that the central part $C(T)$ induces a subtree and contains at least $n-n^{1/3}$ vertices, we can conclude (for sufficiently large $n$) that $v$ lies in $C(T)$. Moreover, since the tree $C(T)$ has no more than 16 leaves, $v$ cannot have more than 16 neighbors in $C(T)$. Those of the components $T_1,T_2,\ldots,T_k$ whose corresponding neighbor of $v$ does not lie in $C(T)$ cannot have more than $n^{1/3}$ vertices in total. Each of the remaining at most 16 components contains at most $\frac{n}{2}$ vertices, so by the pigeonhole principle, there are at least two that contain at least $\frac{n}{32}$ vertices each (for sufficiently large $n$). Without loss of generality, let those be $T_1$ and $T_2$. If the union of all remaining branches (which is $T - (T_1 \cup T_2)$) contains more than $\sqrt{n}$ vertices in total, then we can regard the tree $T$ as the union of $T_1'$, $T_2'$ and $T - (T_1 \cup T_2)$ and apply the argument of Proposition~\ref{prop:centroid} to show that the local mean subtree order at $v$ is at most
$$n - \frac52 \log_2 n + O(1).$$
This gives us a contradiction for sufficiently large $n$. So both $T_1'$ and $T_2'$ contain in fact at least $\frac{n}{2} - \sqrt{n}$ vertices each.

Next, we use~\eqref{eq:s_rec}, which shows that
$$\prod_{i=1}^k s(T_i') = s(T) \leq \sigma(T),$$
regarding $T$ and $T_1',T_2',\ldots$ as rooted at $v$. This implies that either $s(T_1') \leq \sqrt{\sigma(T)}$ or $s(T_2') \leq \sqrt{\sigma(T)}$ (or both). Let us assume that the former holds. We now define $\mathcal{F}_2$ as the family of all subtrees of $T$ that contain $v$, but not all vertices of $T_1 \cap C(T)$. Clearly, $\mathcal{F}_2$ is disjoint from $\mathcal{F}_1$ (whose members need to contain all of $C(T)$).

Note that $T_1 \cap C(T)$ contains at least $|T_1| - n^{1/3} \geq \frac{n}{2} - \sqrt{n} - 1 - n^{1/3} > \frac{n}{3}$ vertices (for sufficiently large $n$). For every $\ell \in \{0,1,\ldots,\lfloor \frac{n}{3} \rfloor \}$, we can find a subtree of $T_1'$ that contains $v$ and $\ell$ vertices of $T_1 \cap C(T)$ (by successively adding vertices). Each of these can be merged with an arbitrary subtree of $T - T_1$ that contains $v$ to form an element of $\mathcal{F}_2$, so we obtain at least
$$\frac{n}{3} \cdot \prod_{i=2}^k s(T_i') = \frac{ns(T)}{3s(T_1')} \geq \frac{n\sigma(T)/(1+n^{-1/4})}{3 \sqrt{\sigma(T)}} = \frac{n\sqrt{\sigma(T)}}{3} \big( 1 - O(n^{-1/4}) \big)$$
such trees, since $s(T) \geq \sigma(T)/(1+n^{-1/4})$ (as we established that $v$ belongs to $C(T)$) and $s(T_1') \leq \sqrt{\sigma(T)}$ by assumption.

Since each such tree does not contain at least $\frac{n}{3} - \ell$ vertices of $T$, we find that the total number of vertices \emph{not} contained in subtrees that belong to $\mathcal{F}_2$ is at least
$$\sum_{\ell=0}^{\lfloor n/3 \rfloor} \Big( \frac{n}{3} - \ell \Big) \cdot \prod_{i=2}^k s(T_i') \geq \frac{n^2\sqrt{\sigma(T)}}{18} \big( 1 - O(n^{-1/4}) \big)$$
by the same inequalities as before.
\end{itemize}
Now we put together the contributions of $\mathcal{F}_1$ and $\mathcal{F}_2$: the total number of vertices not contained in subtrees of $T$ is at least
$$\frac{\sigma(T)}{2} \log_2 \sigma(T) \big(1 - O(n^{-1/4}) \big) + \frac{n^2\sqrt{\sigma(T)}}{18} \big( 1 - O(n^{-1/4}) \big).$$
Thus the mean subtree order of $T$ can be bounded above as follows:
\begin{align*}
\mu_T &\leq n - \frac{1}{\sigma(T)} \Big( \frac{\sigma(T)}{2} \log_2 \sigma(T) + \frac{n^2\sqrt{\sigma(T)}}{18} \Big) \big(1 - O(n^{-1/4}) \big) \\
&\leq n - \frac{1}{2} \log_2 \sigma(T) - \frac{n^2}{20\sqrt{\sigma(T)}} + O(1).
\end{align*}
In the last step, we used the fact that $\log_2 \sigma(T) = O(\log n)$ by Proposition~\ref{prop:subtrees_upper}. Writing $\sigma(T) = x n^4$, we get
$$\mu_T \leq n - 2 \log_2 n - \frac12 \log_2 x - \frac{1}{20\sqrt{x}}.$$
On the other hand,  Theorem~\ref{thm:max_asymp} tells us that $\mu_T = n - 2 \log_2 n + O(1)$. 
As the function $x \mapsto \frac12 \log_2 x + \frac{1}{20\sqrt{x}}$ tends to $\infty$ as $x \to 0$, $x$ must in fact be bounded below by some constant $C_2$, which completes the proof.
\end{proof}

Summarizing, we have shown the following:

\begin{corollary}\label{cor:NumberOfSubtreesInOptimalTree}
The number of subtrees in an optimal tree with $n$ vertices is $\Theta(n^4)$.
\end{corollary}

We remark that the approach that gave us Proposition~\ref{prop:subtrees_upper} and Proposition~\ref{prop:subtrees_lower} could in principle also be used to prove Theorem~\ref{thm:max_asymp}. Of course, the $O$-constants that occur are not nearly optimal.
As the final main result of this section, we are able to provide information on the diameter of optimal trees.

\begin{theorem}\label{thm:diam}
There exists an absolute constant $C$ such that the following statement holds: for every tree $T$ with $n$ vertices and diameter $d$, we have
$$\mu_T \leq n - \log_2 n - 2\log_2 (n-d) + C.$$
\end{theorem}

\begin{proof}
If $d \geq n - \sqrt{n}$, then the statement follows directly from Theorem~\ref{thm:max_asymp}, so we assume that $d \leq n - \sqrt{n}$. Fix a diameter (path of maximum length) of $T$; its length is $d$, so there are $d+1$ vertices on it, and $n-d-1$ vertices that do not lie on it. Let us again consider the central part $C(T)$ as in the previous proof. If $\mu_T <n - 3\log_2 n$, we are done again, so we can assume that $\mu_T \geq n - 3\log_2 n$. Then we can apply the same arguments as in Lemma~\ref{lem:central_part}, where $T$ was assumed to be optimal, but all that was actually used was a lower bound on $\mu_T$. So we can conclude that $C(T)$ contains at least $n - n^{1/3}$ vertices (if $n$ is sufficiently large, as we can always assume). Thus there are at least $n-d - n^{1/3} - 1$ vertices of $C(T)$ that do not lie on the diameter. Each of these lies on at least one path from a leaf of the tree $C(T)$ to the diameter. 

If the tree $C(T)$ has $24$ or more leaves, then we can adapt the argument in the proof of Lemma~\ref{lem:central_part} to show that $\sigma(T) \geq (n^{1/4})^{24} = n^6$, and the proof of Proposition~\ref{prop:subtrees_upper} to show that
$$\mu_T \leq n - \frac12 \log_2 \sigma(T) + O(1) \leq n - 3\log_2 n + O(1) \leq n - \log_2 n - 2\log_2 (n-d) + O(1).$$
In this case we are done, so assume that there are at most $23$ such leaves. By the pigeonhole principle at least one of the paths from a leaf of $C(T)$ to the diameter must have length at least $\frac{1}{23} (n-d - n^{1/3} - 1)$. For sufficiently large $n$, we have $n^{1/3} + 1 \leq \frac12 \sqrt{n} \leq \frac12(n-d)$, so this path has at least length $\frac1{46} (n-d)$.
Let $v$ be the vertex where this path meets the diameter. We know now that there is a path emanating from $v$ that does not have any vertices in common with the diameter other than $v$ and whose length is at least $\frac{1}{46} (n-d)$. Moreover, the vertex $v$ divides the diameter into two pieces. Both pieces must also have a length of at least $\frac{1}{46} (n-d)$, since there would otherwise be a path through $v$ that is longer than the diameter.

We can therefore split $T$ into three subtrees $T_1',T_2',T_3'$ whose union is $T$, whose pairwise intersection is only the vertex $v$, and each of which contains at least $\frac{1}{46} (n-d)$ vertices. The largest of these three subtrees certainly contains at least $\frac{n}{3}$ vertices. We can now use equation~\eqref{eq:mu_add} from the proof of Proposition~\ref{prop:centroid}, which tells us that
$$\mu_T(v) - 1 = \sum_{i=1}^3 \big(\mu_{T_i'}(v) - 1 \big).$$
Now we apply Corollary~\ref{sharplocalresult}:
\begin{align*}
\mu_T(v) - 1 &= \sum_{i=1}^3 \big(\mu_{T_i'}(v) - 1 \big) \\
&\leq \sum_{i=1}^3 \Big( |T_i'| - \log_2 |T_i'| + O(1) \Big) \\
&= |T_1'| + |T_2'| + |T_3'| -  \log_2 |T_1'| -  \log_2 |T_2'| -  \log_2 |T_3'| + O(1) \\
&= |T|   -  \log_2 |T_1'| - \log_2 |T_2'| -  \log_2 |T_3'| + O(1) \\
&\leq n - \log_2 \frac{n}{3} - 2 \log_2 \frac{n-d}{46} + O(1) \\
&= n - \log_2 n - 2\log_2 (n-d) + O(1).
\end{align*}
Note that the $O$-constant does not depend on $d$, so the proof is complete.
\end{proof}

\begin{corollary}\label{cor:diam}
For every positive integer $n$, let $\hat{T}_n$ be an optimal tree with $n$ vertices. Then we have, for every $\delta > 0$,
$$\lim_{n \to \infty} \frac{n - \operatorname{diam} (\hat{T}_n)}{n^{1/2+\delta}} = 0.$$
\end{corollary}

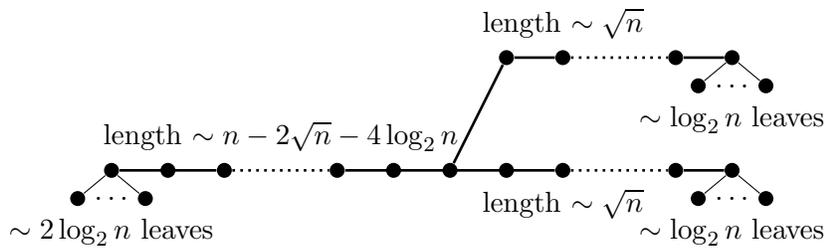
\begin{figure}[htbp]
\centering
    \begin{tikzpicture}[scale=1.5]

        \node[fill=black,circle,inner sep=2pt] (t1) at (0,0) {};
        \node[fill=black,circle,inner sep=2pt] (t2) at (.5,0) {};
        \node[fill=black,circle,inner sep=2pt] (t3) at (1,0) {};

        \node[fill=black,circle,inner sep=2pt] (t5) at (2,0) {};
        \node[fill=black,circle,inner sep=2pt] (t6) at (2.5,0) {};
        \node[fill=black,circle,inner sep=2pt] (t7) at (3,0) {};
        \node[fill=black,circle,inner sep=2pt] (t8) at (3.5,0) {};
        \node[fill=black,circle,inner sep=2pt] (t9) at (4,0) {};
        \node[fill=black,circle,inner sep=2pt] (t11) at (5,0) {};
        \node[fill=black,circle,inner sep=2pt] (t12) at (5.5,0) {};  

        \node[fill=black,circle,inner sep=2pt] (s8) at (3.5,1) {};
        \node[fill=black,circle,inner sep=2pt] (s9) at (4,1) {};
        \node[fill=black,circle,inner sep=2pt] (s11) at (5,1) {};
        \node[fill=black,circle,inner sep=2pt] (s12) at (5.5,1) {};  

\draw [line width=1.pt] (t7)--(s8)--(s9);
\draw [line width=1.pt,dotted] (s9)--(s11);
\draw [line width=1.pt] (s11)--(s12);
               
\draw [line width=1.pt] (t1)--(t3);
\draw [line width=1.pt,dotted] (t3)--(t5);
\draw [line width=1.pt] (t5)--(t9);
\draw [line width=1.pt,dotted] (t9)--(t11);
\draw [line width=1.pt] (t11)--(t12);

    \node[fill=black,circle,inner sep=2pt] (t01) at (-.3,-.25) {};
        \node[fill=black,circle,inner sep=2pt] (t02) at (.3,-.25) {};  
        \node at (0,-.25) {$\ldots$};
\node at (0,-.55) {$\sim 2\log_2 n$ leaves};

\draw (t01)--(t1)--(t02);

    \node[fill=black,circle,inner sep=2pt] (r01) at (5.2,-.25) {};
        \node[fill=black,circle,inner sep=2pt] (r02) at (5.8,-.25) {};  
        \node at (5.5,-.25) {$\ldots$};
\node at (5.5,-.55) {$\sim \log_2 n$ leaves};

\draw (r01)--(t12)--(r02);

    \node[fill=black,circle,inner sep=2pt] (s01) at (5.2,.75) {};
        \node[fill=black,circle,inner sep=2pt] (s02) at (5.8,.75) {};  
        \node at (5.5,.75) {$\ldots$};
\node at (5.5,.45) {$\sim \log_2 n$ leaves};

\draw (s01)--(s12)--(s02);

\node at (4,-.3) {length $\sim \sqrt{n}$};
\node at (4,1.3) {length $\sim \sqrt{n}$};

\node at (1.5,.3) {length $\sim n - 2\sqrt{n} - 4 \log_2 n $};

        \end{tikzpicture}
\caption{Example of a tree with near-maximum mean subtree order.}\label{fig:maxden}
\end{figure}

In plain words, optimal trees must have a diameter that is close to the number of vertices. We remark that there are trees whose diameter is about $n - \sqrt{n}$ for which the asymptotic formula of Theorem~\ref{thm:max_asymp} is attained. The construction is fairly simple: merge three brooms at their roots; two of them have length approximately $\sqrt{n}$ and $\log_2 n$ leaves each. The third one consists of a path of length about $n - 2\sqrt{n} - 4\log_2 n$, with approximately $2 \log_2 n$ leaves at the end (Figure~\ref{fig:maxden}). It is not difficult to check that the resulting tree satisfies the abovementioned properties: the diameter is $n - \sqrt{n} + O(\log n)$, the mean subtree order is $n - 2\log_2 n + O(1)$.

\section{Constructing better caterpillars}\label{sec:opt_cat}


In this section, we analyze a construction of caterpillars that achieve a slightly higher mean subtree order than double-brooms.  This allows us to improve the final constant in the lower bound for the maximum mean subtree order.
We conjecture that our construction in fact accurately reflects the shape of optimal trees for large $n$.

For $n \le 24$ the optimal trees have been computed explicitly. In particular, see \cite[Figure~1]{mol2019} for the optimal trees when $16 \le n \le 24$. None of these are double-brooms.

\begin{theorem}\label{thr:optimaltreestructure}
	For large $n$, there is always a tree $T$ with $n$ vertices such that $\mu_{T} \geq n - 2\log_2 (0.9n) + f(2\log_2 (0.9n) ) + o(1)$, where $f$ is the same $1$-periodic function as in Corollary~\ref{sharplocalresult}, given by $f(x) = x - 2^x$ for $x \in [0,1]$.
\end{theorem}

\begin{proof}
We provide an explicit construction.	Let $T$ be a caterpillar consisting of a path with $\ell+1$ vertices (which will be called the \emph{stem}), $m$ leaves attached at either end, and $k$ additional leaves, where $n=\ell+2m+k+1$. We will choose $m$ and $k$ in such a way that the number of leaves is $2m + k = 4 \log_2 n + O(1)$ and $k = c \log_2 n + O(1)$ for some fixed constant $c \in (0,1)$.
	
Fix an end of the stem and call it the left end; the other one will be called the right end. The $k$ additional leaves are attached to vertices of the stem whose distances from the left end are
$a_1 \ell, a_2 \ell, \ldots,  a_k \ell$ respectively, where $0 \leq a_1 \leq a_2 \leq \ldots \leq a_k \leq 1$. These will be called \emph{support vertices}. Moreover, we set $a_0=0$ and $a_{k+1}=1$.

Let us now determine the number of subtrees as well as their total order. Firstly, there are $2^{2m+k}$ subtrees that contain the entire stem. Next, we consider subtrees containing the left end, but not the right. Here, the number of subtrees containing precisely the first $i$ support vertices is
$$2^{m+i} (a_{i+1}\ell - a_i \ell),$$
since there are $m+i$ leaves that can potentially be added, and $a_{i+1}\ell - a_i \ell$ ways to add a path segment between the $i$-th and the $(i+1)$-th support vertex. Using the same reasoning, we find that there are 
$$2^{m+k-i} (a_{i+1}\ell - a_i \ell)$$
subtrees containing the right end, but not the left, and precisely the last $k-i$ support vertices.

Finally, the number of subtrees containing neither of the two ends of the stem can be bounded by $O(2^k \ell^2)$, as those are either single leaves or consist of part of the stem and a subset of the $k$ additional leaves.
So the total number of subtrees of $T$ is 
	\begin{align*}
	\sigma(T)&=2^{2m+k}+2^m\ell \left( \sum_{i=0}^{k} (2^i+2^{k-i}) (a_{i+1}-a_i) \right) + O(2^k \ell^2)\\
	&=2^{2m+k}+2^m\ell \left( 1+2^k +\sum_{i=1}^{k} (2^{k-i}-2^{i-1}) a_i   \right) + O(2^k \ell^2).
	\end{align*}
	
In a similar fashion, we can compute $\tau(T)$, the sum of the orders of the subtrees of $T$. The subtrees containing both ends of the stem contribute $2^{2m+k}(\ell+1+m+\frac k2) = 2^{2m+k} (n - m - \frac k2)$. The subtrees that contain the left end and the first $i$ support vertices, but not the right end, contribute a total of
$$2^{m+i} (a_{i+1}\ell - a_i \ell) \frac{m+i+1 + a_i \ell + a_{i+1} \ell}{2} = 2^{m+i-1} \ell^2 (a_{i+1}^2 - a_i^2) + O\Big((m+k) 2^{m+k} (a_{i+1}-a_i) \ell \Big).$$
Likewise, subtrees that contain the right end and the last $k-i$ support vertices, but not the left end, contribute a total of
\begin{multline*}
2^{m+k-i} (a_{i+1}\ell - a_i \ell) \frac{m+k-i+1 + (1-a_i)\ell +(1- a_{i+1}) \ell}{2} \\
= 2^{m+k-i-1} \ell^2 (2a_{i+1}-2a_i + a_i^2 - a_{i+1}^2) + O\Big((m+k) 2^{m+k} (a_{i+1}-a_i) \ell \Big).
\end{multline*}
The contribution of all subtrees that do not contain either of the ends is bounded by $O(2^k \ell^3)$ by the same reasoning as before. So we have 
\begin{align*}
\tau(T) &= 2^{2m+k}\Big(n - m - \frac{k}{2} \Big)+2^m\ell^2 \left( \sum_{i=0}^{k} (2^{i-1} - 2^{k-i-1}) (a^2_{i+1}-a^2_i)+\sum_{i=0}^{k} 2^{k-i} (a_{i+1}-a_i)\right) \\
&\quad+ O\Big((m+k)2^{m+k}\ell + 2^k \ell^3\Big) \\
	&=2^{2m+k}\Big(n - m - \frac{k}{2} \Big)+2^m\ell^2 \left( 2^{k-1}+\frac12 - \sum_{i=1}^{k} (2^{k-i-1}+2^{i-2}) a_i^2 +\sum_{i=1}^{k} 2^{k-i} a_i \right) \\
&\quad+ O\Big((m+k)2^{m+k}\ell + 2^k \ell^3\Big).
	\end{align*}
Recall now the assumptions we made on $m$ and $k$. With those, we find that $2^{2m+k} = \Theta(n^4)$ and $2^k \ell^2 = \Theta(n^{2+c})$ as well as $\ell = n - O(\log n)$. Consequently,
\begin{align*}
\sigma(T) &= 2^{2m+k} + 2^{m+k} \ell p(k;a_1,a_2,\ldots,a_k) + O(2^k \ell^2) \\
&= 2^{2m+k} \Big(1 + 2^{-m} n p(k;a_1,a_2,\ldots,a_k) + O(n^{c-2}) \Big),
\end{align*}
where
$$p(k;a_1,a_2,\ldots,a_k) = 2^{-k}+1 +\sum_{i=1}^{k} (2^{-i}-2^{i-k-1}) a_i $$
is easily seen to be bounded. Similarly,
$$\tau(T) = 2^{2m+k}\Big(n-m-\frac k2 + 2^{-m} n^2 q(k;a_1,a_2,\ldots,a_k) + O(n^{c-1})\Big),$$
where
$$q(k;a_1,a_2,\ldots,a_k) = \frac12 + 2^{-k-1} - \sum_{i=1}^{k} (2^{-i-1}+2^{i-k-2}) a_i^2 +\sum_{i=1}^{k} 2^{-i} a_i$$
is bounded. Thus we obtain
$$\mu_T = \frac{\tau(T)}{\sigma(T)} = n-m-\frac k2 + 2^{-m} n^2 \big( q(k;a_1,a_2,\ldots,a_k) - p(k;a_1,a_2,\ldots,a_k) \big) + O(n^{c-1}).$$
We can still choose $a_1,a_2,\ldots,a_k$, and we make the choice in such a way that $q(k;a_1,a_2,\ldots,a_k) - p(k;a_1,a_2,\ldots,a_k)$ is (near) maximal. To this end, note that
$$ q(k;a_1,a_2,\ldots,a_k) - p(k;a_1,a_2,\ldots,a_k) = 
- \frac12 - 2^{-k-1} - \sum_{i=1}^{k} \Big( (2^{-i-1}+2^{i-k-2}) a_i^2 -2^{i-k-1} a_i \Big).$$
The sum of quadratic polynomials attains its maximum when $a_i= \frac{1}{2^{k-2i+1}+1}$ for every $i$. However, since $a_i \ell$ needs to be an integer for each $i$, we can only choose it in such a way that $a_i = \frac{1}{2^{k-2i+1}+1} + O(n^{-1})$. The error term here has an asymptotically negligible impact on $\mu_T$. With this choice, we arrive at
$$\mu_T = n-m-\frac k2 - 2^{-m} n^2 \Big( \frac12 + 2^{-k-1} - \sum_{i=1}^{k} \frac{2^{i-k-2}}{2^{k-2i+1}+1} \Big) + O(n^{c-1}).$$
The expression in brackets simplifies as follows:
\begin{align*}
\frac12 + 2^{-k-1} - \sum_{i=1}^{k} \frac{2^{i-k-2}}{2^{k-2i+1}+1}
&= \frac12 + 2^{-k-1} - \sum_{i=1}^{k} 2^{i-k-2} + \sum_{i=1}^{k} \frac{2^{-i-1}}{2^{k-2i+1}+1} \\
&= 2^{-k} + \sum_{i=1}^{k} \frac{2^{-i-1}}{2^{k-2i+1}+1} \\
&= 2^{-k} +  \sum_{i=1}^{k} \frac{2^{-(k+3)/2}}{2^{(k+1)/2-i}+2^{i-(k+1)/2}} \\
&= 2^{-k/2} \sum_{i=-\infty}^{\infty}  \frac{2^{-3/2}}{2^{(k+1)/2-i}+2^{i-(k+1)/2}} + O(2^{-k}).
\end{align*}
The value of the infinite sum depends only on the residue class of $k$ modulo $2$: it is approximately $0.801214$ for even $k$
and approximately $0.801218$ for odd $k$. In particular, it is less than $0.81$. Therefore, we have
$$\mu_T \geq n - m - \frac k2 - 0.81 \cdot 2^{-m-k/2} n^2 + O(n^{c-1}) = n - m' - 0.81 \cdot 2^{-m'} n^2 + O(n^{c-1}),$$
where $m' = m + \frac k2$. It remains to minimize the expression
$$m' + 0.81 \cdot 2^{-m'} n^2 = m' + \frac{(0.9n)^2}{2^{m'}},$$
and (assuming for simplicity that $k$ is chosen to be even) we already know from the proof of Corollary~\ref{sharplocalresult} that this can be achieved by taking $m' = \lfloor 2\log_2 (0.9n) \rfloor$, resulting in the lower bound
$$\mu_T \geq n - 2\log_2 (0.9n)+ f(2\log_2 (0.9n)) + o(1),$$
which completes the proof.
\end{proof}

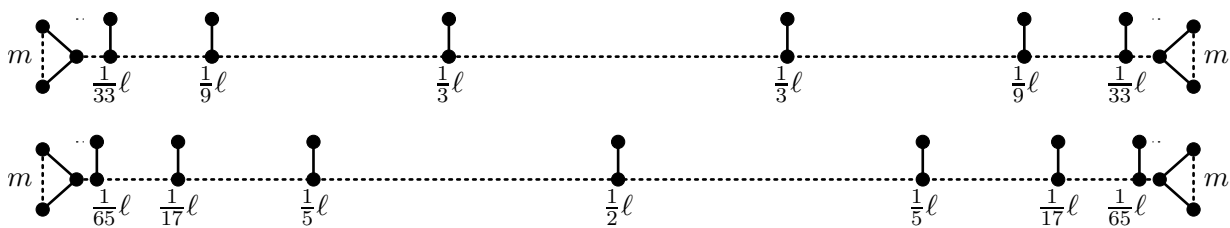
\begin{figure}[h]
	\centering

	\begin{tikzpicture}[line cap=round,line join=round,>=triangle 45,x=0.9cm,y=1.0cm]
	\clip(-0.5,-0.8) rectangle (17.7,0.8);

	\draw [line width=0.5pt,dotted] (0.5,0.5)--(0.6,0.5);
	\draw [line width=0.5pt,dotted] (16.5,0.5)--(16.4,0.5);
	
\draw [line width=1.pt,dotted] (0.5,0)--(16.5,0);

\draw [line width=1.pt] (1,0)--(1,0.5);
\draw [line width=1.pt] (2.5,0)--(2.5,0.5);

\draw [line width=1.pt] (6,0)--(6,0.5);
\draw [line width=1.pt] (11,0)--(11,0.5);
	
\draw [line width=1.pt] (16,0)--(16,0.5);
\draw [line width=1.pt] (14.5,0)--(14.5,0.5);
	
	\draw [line width=1.pt] (0,0.4)--(0.5,0)-- (0,-0.4);
	\draw [line width=1.pt] (17,0.4)--(16.5,0)-- (17,-0.4);
	\draw [line width=1.pt,dotted] (0,0.4)-- (0,-0.4);
	\draw [line width=1.pt,dotted] (17,0.4)-- (17,-0.4);
	\draw (0,0) node[anchor= east] {$m$};		
	\draw (17,0) node[anchor=west] {$m$};

	\draw (1,0) node[anchor=north] {$\frac 1{33} \ell$};
\draw (2.5,0) node[anchor=north] {$\frac 19 \ell$};
\draw (6,0) node[anchor=north] {$\frac 13 \ell$};
	\draw (16,0) node[anchor=north] {$\frac 1{33} \ell$};
\draw (14.5,0) node[anchor=north] {$\frac 19 \ell$};
\draw (11,0) node[anchor=north] {$\frac 13 \ell$};

	\begin{scriptsize}
	\draw [fill=black] (0,0.4) circle (2.5pt);
	\draw [fill=black] (0,-0.4) circle (2.5pt);
	\draw [fill=black] (17,0.4) circle (2.5pt);
	\draw [fill=black] (17,-0.4) circle (2.5pt);
	
	\draw [fill=black] (0.5,0) circle (2.5pt);
	
	\draw [fill=black] (2.5,0) circle (2.5pt);
	\draw [fill=black] (1,0) circle (2.5pt);	
	\draw [fill=black] (6,0) circle (2.5pt);
	\draw [fill=black] (11,0) circle (2.5pt);
	\draw [fill=black] (14.5,0) circle (2.5pt);
	\draw [fill=black] (16,0) circle (2.5pt);
	
	\draw [fill=black] (2.5,0.5) circle (2.5pt);
	\draw [fill=black] (1,0.5) circle (2.5pt);	
	\draw [fill=black] (6,0.5) circle (2.5pt);
	\draw [fill=black] (11,0.5) circle (2.5pt);
	\draw [fill=black] (14.5,0.5) circle (2.5pt);
	\draw [fill=black] (16,0.5) circle (2.5pt);

	
	\draw [fill=black] (16.5,0) circle (2.5pt);
	
	\end{scriptsize}
	\end{tikzpicture}
	\quad
	\begin{tikzpicture}[line cap=round,line join=round,>=triangle 45,x=0.9cm,y=1.0cm]
\clip(-0.5,-0.8) rectangle (17.7,0.8);

\draw [line width=1.pt,dotted] (0.5,0)--(16.5,0);

\draw [line width=0.5pt,dotted] (0.5,0.5)--(0.6,0.5);
\draw [line width=0.5pt,dotted] (16.5,0.5)--(16.4,0.5);

\draw [line width=1.pt] (8.5,0)--(8.5,0.5);

\draw [line width=1.pt] (0.8,0)--(0.8,0.5);
\draw [line width=1.pt] (2,0)--(2,0.5);

\draw [line width=1.pt] (4,0)--(4,0.5);
\draw [line width=1.pt] (13,0)--(13,0.5);

\draw [line width=1.pt] (16.2,0)--(16.2,0.5);
\draw [line width=1.pt] (15,0)--(15,0.5);

\draw [line width=1.pt] (0,0.4)--(0.5,0)-- (0,-0.4);
\draw [line width=1.pt] (17,0.4)--(16.5,0)-- (17,-0.4);
\draw [line width=1.pt,dotted] (0,0.4)-- (0,-0.4);
\draw [line width=1.pt,dotted] (17,0.4)-- (17,-0.4);
\draw (0,0) node[anchor= east] {$m$};		
\draw (17,0) node[anchor=west] {$m$};

\draw (8.5,0) node[anchor=north] {$\frac1{2}  \ell$};
\draw (4,0) node[anchor=north] {$\frac1{5}  \ell$};
\draw (2,0) node[anchor=north] {$\frac1{17}  \ell$};
\draw (1,0) node[anchor=north] {$\frac1{65}  \ell$};

\draw (13,0) node[anchor=north] {$\frac1{5}  \ell$};
\draw (15,0) node[anchor=north] {$\frac1{17}  \ell$};
\draw (16,0) node[anchor=north] {$\frac1{65}  \ell$};

\begin{scriptsize}
	\draw [fill=black] (0,0.4) circle (2.5pt);
	\draw [fill=black] (0,-0.4) circle (2.5pt);
	\draw [fill=black] (17,0.4) circle (2.5pt);
	\draw [fill=black] (17,-0.4) circle (2.5pt);
	
	\draw [fill=black] (0.5,0) circle (2.5pt);
	
	\draw [fill=black] (8.5,0) circle (2.5pt);
	\draw [fill=black] (8.5,0.5) circle (2.5pt);
	
	\draw [fill=black] (2,0) circle (2.5pt);
	\draw [fill=black] (0.8,0) circle (2.5pt);	
	\draw [fill=black] (4,0) circle (2.5pt);
	\draw [fill=black] (13,0) circle (2.5pt);
	\draw [fill=black] (15,0) circle (2.5pt);
	\draw [fill=black] (16.2,0) circle (2.5pt);
	
	\draw [fill=black] (2,0.5) circle (2.5pt);
	\draw [fill=black] (0.8,0.5) circle (2.5pt);	
	\draw [fill=black] (4,0.5) circle (2.5pt);
	\draw [fill=black] (13,0.5) circle (2.5pt);
	\draw [fill=black] (15,0.5) circle (2.5pt);
	\draw [fill=black] (16.2,0.5) circle (2.5pt);

	
	\draw [fill=black] (16.5,0) circle (2.5pt);
	
\end{scriptsize}
\end{tikzpicture}
	\caption{A sketch of the trees constructed in the proof of Theorem~\ref{thr:optimaltreestructure}: even $k$ (top) and odd $k$ (bottom).}
	\label{fig:limittrees?}
\end{figure}

Figure~\ref{fig:limittrees?} shows the rough structure of the trees that are constructed in the proof of Theorem~\ref{thr:optimaltreestructure}. Note the wide variety of potential choices for both $k$ and $a_1,a_2,\ldots,a_k$, which results in a large number of trees that satisfy the asymptotic inequality in  Theorem~\ref{thr:optimaltreestructure}.

As an immediate corollary, we obtain the following result.

\begin{corollary}
	For sufficiently large $n$, no double-broom is an optimal tree.
\end{corollary}

\begin{proof}
The difference in mean subtree order between the trees constructed in Theorem~\ref{thr:optimaltreestructure} and the best double-broom is at least 
$$-2\log_2(0.9) + f(2\log_2(0.9n)) - f(2\log_2 n) + o(1),$$
which can be expressed as $g(2\log_2 n) + o(1)$ for a $1$-periodic function given by
$$g(x) = \max \big( 0.19 \cdot 2^x, 1 - 0.62 \cdot 2^x \big)$$
for $x \in [0,1]$. Since the minimum of this function is positive (approximately $0.234568$), the statement follows.
\end{proof}

By means of a computer program\footnote{\url{https://github.com/StijnCambie/Jamison}, document Check\_NoDoubleBroom.}, it can also be checked that for $25 \leq n \leq 1000$, the best balanced double-broom is not optimal either.

\section{Optimal trees are nearly caterpillars}\label{sec:NearlyCatterpillar}

In this concluding section, we summarize the progress towards the last major open question by Jamison on the mean subtree order of trees. Although there is not yet a proof for the optimal tree to be a caterpillar, the evidence listed below shows that the optimal trees are indeed very much like caterpillars.
In Corollary~\ref{cor:NumberOfSubtreesInOptimalTree}, we proved that the number of subtrees in an optimal tree is of the same order as the number of subtrees in the optimal double broom. 
To further analyze the structure of the optimal tree $\hat{T}_n$, we combine Corollary~\ref{cor:boundsmaxmu} and Theorem~\ref{thm:diam} to see that 
$$ n- \operatorname{diam}( \hat{T}_n) = O( \sqrt n). $$
In particular, this implies that almost all, up to at most $O( \sqrt n)$, vertices belong to one path.
We conclude with a final observation. Set $d = \operatorname{diam}( \hat{T}_n)$, let $e$ denote an end vertex of the diameter, and fix some constant $0<\eps <\frac 12$. Let $v$ be a vertex on the diameter with $\eps d< d(e,v)<(1-\eps)d$, and let $T'$ be a subtree of $\hat{T}_n$ that only shares the vertex $v$ with the diameter. If $\hat{T}_n$ is a caterpillar, then $T'$ has to be a star. We prove the weaker statement that $T'$ cannot contain more than $C_{\eps}$ vertices, where $C_{\eps}$ is a constant of the order $\frac{1}{\eps(1-\eps)}$. To see this, note that $\hat{T}_n$ can be divided into three subtrees $T_1$, $T_2$ and $T'$ that share only the vertex $v$, where $|T_1|, |T_2| \geq \eps d$ and $|T_1| + |T_2| \geq d$. It follows that the sum of the defects is (by Theorem~\ref{sharplocalresult}) at least
\begin{align*}
\log_2 |T_1| + \log_2 |T_2| + \log_2 |T'| + O(1) &\geq 2 \log_2 d + \log_2 (\eps(1-\eps)) + \log_2 |T'| + O(1) \\
&= 2 \log_2 n + \log_2 (\eps(1-\eps)) + \log_2 |T'| + O(1).
\end{align*}
This has to be smaller than $2 \log_2 n -f(2 \log_2 n ) + o(1)$ by Theorem~\ref{thm:max_asymp}, so it follows that $\log_2 |T'| \leq \log_2 \frac{1}{\eps(1-\eps)} + O(1)$.

\subsection*{Acknowledgment}
This work originated from a minisymposium at CanaDAM 2019, to which the authors were invited by Lucas Mol and Ortrud Oellermann. We would like to thank them for inviting us.
The first author was supported by a grant with grant number 00352039192000OO for attending CanaDAM.
The authors are also very grateful to the anonymous referees for the time they dedicated to this article and for the helpful comments they made.

\bibliographystyle{abbrv}
\bibliography{MeanSubtreeOrder}

\end{document}